\documentclass[a4paper,12pt]{amsart}

\usepackage{float}
\usepackage{euscript,eufrak,verbatim}
\usepackage{graphicx}
\usepackage[usenames]{color}
\usepackage[colorlinks,linkcolor=red,anchorcolor=blue,citecolor=blue]{hyperref}
\usepackage{amsmath}
\usepackage{amsthm}
\usepackage[all]{xy}
\usepackage{graphicx} 
\usepackage{amssymb}

\theoremstyle{plain}
\newtheorem{main theorem}{Main Theorem}
\newtheorem{theorem}{Theorem}[section]
\newtheorem{lemma}[theorem]{Lemma}

\newtheorem{corollary}[theorem]{Corollary}
\newtheorem{proposition}[theorem]{Proposition}
\newtheorem{claim}[theorem]{Claim}

\theoremstyle{definition}
\newtheorem{definition}[theorem]{Definition}
\newtheorem{remark}[theorem]{Remark}

\newtheorem{problem}[theorem]{Problem}

\makeatletter 
\@addtoreset{equation}{section}
\numberwithin{equation}{section}

\newcommand{\norm}[1]{\left\lVert#1\right\rVert}

\newcommand{\diam}{\mathrm{Diam}}
\newcommand{\supp}{\mathrm{supp}}
\newcommand{\mdim}{\mathrm{mdim}}

\newcommand{\widim}{\mathrm{Widim}}
\newcommand{\st}{\mathrm{St}}
\newcommand{\bc}{\mathrm{bc}}

\newcommand{\ocap}{\mathrm{ocap}}

\title{On an analogue of the Hurewicz theorem for mean dimension}

\author{Masaki Tsukamoto}

\address
{Department of Mathematics, Kyoto University, Kitashirakawa Oiwake-cho, Sakyo-ku, Kyoto 606-8502, Japan}

\email{tsukamoto@math.kyoto-u.ac.jp}

\begin{document}



\thanks{M.T. was supported by JSPS KAKENHI JP21K03227.}

\maketitle

\begin{abstract}
The Hurewicz theorem is a fundamental result in classical dimension theory concerning 
continuous maps which lower topological dimension.
We study whether or not its analogue holds for mean dimension of dynamical systems.
Our first main result shows that an analogue of the Hurewicz theorem does not hold for mean dimension in 
general.
Our second main result shows that it holds true if a base system has zero mean dimension.
\end{abstract}

\section{Introduction} \label{section: introduction}

\subsection{Background and basic definitions}  \label{subsection: background and basic definitions}

The Hurewicz theorem is one of the fundamental results in the classical topological dimension theory.
For a compact metrizable space $X$ we denote its topological dimension (Lebesgue covering dimension) by $\dim X$.
Let $f:X\to Y$ be a continuous map between compact metrizable spaces. 
Then the Hurewicz theorem \cite[p. 91, Theorem VI 7]{Hurewicz--Wallman} states that 
\[ \dim X \leq \dim Y + \sup_{y\in Y} \dim f^{-1}(y). \]
The purpose of this paper is to study an analogue of this theorem for mean dimension of dynamical systems.

Mean dimension is a dynamical version of topological dimension introduced by Gromov \cite{Gromov}.
It quantifies how many parameters per iterate we need for describing orbits of a dynamical system.
It has applications to several problems in topological dynamics \cite{Lindenstrauss--Weiss, Lindenstrauss, Meyerovitch--Tsukamoto, Gutman--Tsukamoto minimal}.

We need to prepare some definitions before rigorously stating our problem.
Throughout of this paper we assume that “simplicial complex” 
means a \textit{finite} simplicial complex (namely, the number of its simplices is finite).
Let $(X, d)$ be a metric space and $Y$ a topological space. 
A continuous map $f:X\to Y$ is called an \textbf{$\varepsilon$-embedding} if we have $\diam f^{-1}(y) < \varepsilon$ for all $y\in Y$.
We define the \textbf{$\varepsilon$-width dimension}
$\widim_\varepsilon(X, d)$ as the minimum integer $n\geq 0$ such that there exist an $n$-dimensional simplicial complex $P$ and
an $\varepsilon$-embedding $f:X\to P$.
If $X$ is compact then its topological dimension is defined by 
\[  \dim X  = \lim_{\varepsilon\to 0} \widim_\varepsilon(X, d). \]

A pair $(X, T)$ is called a \textbf{dynamical system} if $X$ is a compact metrizable space and $T:X\to X$ is a homeomorphism.
Let $(X,T)$ be a dynamical system with a metric $d$ on $X$.
For each natural number $N$ we define a metric $d_N$ on $X$ by 
\[  d_N(x, y) = \max_{0\leq n<N} d\left(T^n x, T^n y\right).  \]
We define the \textbf{mean dimension} of $(X, T)$ by 
\[  \mdim(X, T) = \lim_{\varepsilon\to 0}\left(\lim_{N\to \infty} \frac{\widim_\varepsilon(X, d_N)}{N}\right). \]
The value of $\mdim(X, T)$ is independent of the choice of a metric $d$, and it provides a topological invariant of $(X,T)$.

Let $A$ be a (not necessarily invariant) closed subset of $X$. 
We define the \textbf{upper and lower mean dimensions} of $A$ by 
\begin{align*}
     \overline{\mdim}(A,T) &
     = \lim_{\varepsilon\to 0} \left(\limsup_{N\to \infty} \frac{\widim_\varepsilon(A,d_N)}{N}\right), \\
      \underline{\mdim}(A,T) &
      = \lim_{\varepsilon\to 0} \left(\liminf_{N\to \infty} \frac{\widim_\varepsilon(A,d_N)}{N}\right). 
\end{align*}      
These are also independent of the choice of $d$.

Let $(X, T)$ and $(Y, S)$ be dynamical systems.
A map $\pi:X\to Y$ is called a \textbf{factor map} between dynamical systems if $\pi$ is a continuous surjection satisfying 
$\pi\circ T = S\circ \pi$. 
We often denote it by $\pi:(X, T)\to (Y, S)$ for clarifying the underlying dynamics.

Let $\pi:(X, T)\to (Y, S)$ be a factor map between dynamical systems.
We would like to study the mean dimension of fibers of $\pi$.
For this purpose we define the \textbf{relative mean dimension} of $\pi:(X, T)\to (Y, S)$ by 
\[ \mdim(\pi, T) = \lim_{\varepsilon\to 0} 
      \left(\lim_{N\to \infty} \frac{\sup_{y\in Y} \widim_\varepsilon\left(\pi^{-1}(y), d_N\right)}{N}\right). \]
It is easy to check that the quantity $\sup_{y\in Y} \widim_\varepsilon\left(\pi^{-1}(y), d_N\right)$ is sub-additive in $N$ 
and monotone in $\varepsilon$. So the above limits exist.

The next proposition clarifies the meaning of this definition.

\begin{proposition}  \label{prop: relative mean dimension}
Let $\pi:(X, T)\to (Y, S)$ be a factor map between dynamical systems. Then
\[  \mdim(\pi, T) = \sup_{y\in Y} \overline{\mdim}\left(\pi^{-1}(y), T\right) = \sup_{y\in Y} \underline{\mdim}\left(\pi^{-1}(y), T\right). \]
\end{proposition}
So we can say that the relative mean dimension $\mdim(\pi, T)$ properly measures the mean dimension of fibers of $\pi$.

Now we formally state the main problem we study in the paper:

\begin{problem}  \label{main problem}
Let $\pi:(X, T)\to (Y, S)$ be a factor map between dynamical systems. 
Does the following inequality hold true?
\begin{equation} \label{eq: Hurewicz inequality for mean dimension}
      \mdim(X, T) \leq \mdim(Y, S) + \mdim(\pi, T). 
\end{equation}      
\end{problem}
This problem was originally posed by the author in \cite[Problem 4.8]{Tsukamoto moduli space of Brody curves}
more than ten years ago.
He encountered it when he studied mean dimension of certain dynamical systems coming from 
geometric analysis.
At that time he spent a lot of time trying to prove the inequality (\ref{eq: Hurewicz inequality for mean dimension}),
but he did not succeed.
He could only prove an inequality much weaker than (\ref{eq: Hurewicz inequality for mean dimension}) under a rather 
artificial assumption \cite[Theorem 4.6]{Tsukamoto moduli space of Brody curves}.

Recently Liang \cite{Liang} revisited this problem from a new angle.
He proved, among other things, that the inequality (\ref{eq: Hurewicz inequality for mean dimension}) holds true in the category of 
\textit{algebraic actions}\footnote{Indeed, in the category of algebraic actions, 
the inequality (\ref{eq: Hurewicz inequality for mean dimension}) becomes an equality.
This is also proved in \cite[Corollary 2.18]{Liang}.}
\cite[Corollary 2.18]{Liang}.
See also \cite[Corollary 6.1]{Li--Liang} for a closely related result.

We will give both (partially) positive and negative answers to Problem \ref{main problem} below.

\subsection{Main results}  \label{subsection: main results}

Our first main result shows a negative answer to Problem \ref{main problem} in a rather strong sense:

\begin{theorem}  \label{theorem: Hurewicz does not hold in general}
For any positive number $\delta$ there exists a factor map between dynamical systems
\[  \pi: (X, T)\to (Y, S) \]
satisfying
\[  \mdim(X, T) = 1, \quad  \mdim(Y, S) < \delta,  \quad \mdim(\pi, T) = 0. \] 
\end{theorem}

In particular, letting $\delta<1$, this shows that the inequality 
\[  \mdim(X, T) \leq \mdim(Y, S) + \mdim(\pi, T) \]
does not hold in general.

\begin{remark}  \label{remark: raising mean dimension of (X,T)}
In the above statement, we consider the condition $\mdim(X, T) = 1$.
This is just for simplicity, and indeed
we can make $\mdim(X, T)$ arbitrary large as follows:
Let $n$ be a natural number and $\delta$ a positive number.
By Theorem \ref{theorem: Hurewicz does not hold in general} there exists a factor map $\pi:(X, T)\to (Y, S)$ satisfying 
\[ \mdim(X, T) = 1, \quad  \mdim(Y, S) < \frac{\delta}{n},  \quad \mdim(\pi, T) = 0. \] 
Then it is easy to check that the factor map $\pi:(X, T^n) \to (Y, S^n)$ satisfies 
\[  \mdim(X, T^n) = n, \quad \mdim(Y, S^n) < \delta, \quad \mdim(\pi, T^n) = 0. \]
Furthermore we can even prove that for any positive number $\delta$ there exists a
factor map $\pi^\prime:(X^\prime, T^\prime)\to (Y^\prime, S^\prime)$ satisfying 
\[  \mdim(X^\prime, T^\prime) = \infty, \quad  \mdim(Y^\prime, S^\prime) < \delta,  \quad \mdim(\pi^\prime, T^\prime) = 0. \] 
See Remarks \ref{remark: raising mean dimension} and \ref{remark: universality} 
in \S \ref{subsection: Proof of Theorem Hurewicz does not hold} for further discussions.
\end{remark}

Some readers might wonder whether one can even require $\mdim(Y, S) = 0$ instead of $\mdim(Y, S) < \delta$
in the statement of Theorem \ref{theorem: Hurewicz does not hold in general}.
However this turns out to be impossible. This is our second main result:

\begin{theorem}  \label{theorem: Hurewicz holds if the base is zero mean dimensional}
Let $\pi:(X, T)\to (Y, S)$ be a factor map between dynamical systems.
If $\mdim(Y, S)=0$ then
\[  \mdim(X, T) = \mdim(\pi, T). \]
\end{theorem}

Namely, an analogue of the Hurewicz theorem holds true if the base system $(Y, S)$ has zero mean dimension.
This provides a partially positive answer to Problem \ref{main problem}.

The proof of Theorem \ref{theorem: Hurewicz does not hold in general} is based on a result of 
Gromov \cite{Gromov width} and its variations developed in \S \ref{subsection: variations of Gromov’s lemma} below.
A main ingredient of the proof of Theorem \ref{theorem: Hurewicz holds if the base is zero mean dimensional} is 
Lindenstrauss--Weiss’ theory of small boundary property.

\textbf{Acknowledgment.}
I would like to deeply thank Tom Meyerovitch.
The conversation with him revived my interest in Problem \ref{main problem}.
We discussed the problem together.
He generously declined to be a coauthor of this paper.
But his influence is visible in several places of the paper.
In particular I came up with Remark \ref{remark: universality} thanks to his suggestion.

\section{Proof of Proposition \ref{prop: relative mean dimension}}   \label{section: proof of Proposition relative mean dimension}

\subsection{Preparations on $\varepsilon$-width dimension} \label{subsection: preparations on widim}

Here we prepare some simple results on $\varepsilon$-width dimension.
Let $(X, d)$ and $(Y, d^\prime)$ be metric spaces. We consider its product 
\[  (X, d) \times (Y, d^\prime)  = \left(X\times Y, d\times d^\prime\right), \]
where $d\times d^\prime$ is a metric on $X\times Y$ defined by 
\[  d\times d^\prime\left((x_1,y_1), (x_2, y_2)\right) = \max\left(d(x_1, x_2), d^\prime(y_1, y_2)\right). \]

\begin{lemma}  \label{lemma: product and widim}
For any $\varepsilon>0$ we have 
\[  \widim_\varepsilon\left(X\times Y, d\times d^\prime\right) \leq 
     \widim_\varepsilon(X, d) + \widim_\varepsilon(Y, d^\prime). \]
\end{lemma}

\begin{proof}
If $f:(X, d)\to K$ and $g:(Y, d^\prime) \to L$ are both $\varepsilon$-embeddings then
$f\times g: (X\times Y, d\times d^\prime) \to K\times L$ is also an $\varepsilon$-embedding.
\end{proof}

We say that a map $f:X\to Y$ is \textbf{distance non-decreasing} if for every $x_1, x_2\in X$ we have 
\[  d(x_1, x_2) \leq d^\prime\left(f(x_1), f(x_2)\right). \] 

\begin{lemma} \label{lemma: distance non-decreasing map and widim}
If there exists a distance non-decreasing continuous map $f:X\to Y$ then for any $\varepsilon>0$
\[  \widim_\varepsilon\left(X, d\right) \leq \widim_\varepsilon(Y, d^\prime). \]
\end{lemma}

\begin{proof}
If $g:Y\to K$ is an $\varepsilon$-embedding then $g\circ f:X\to K$ is also.
\end{proof}

\subsection{Proof of Proposition \ref{prop: relative mean dimension}} 
\label{subsection: proof of Proposition relative mean dimension}
Let $\pi:(X, T)\to (Y, S)$ be a factor map between dynamical systems.
Here we prove Proposition \ref{prop: relative mean dimension}.
It is obvious from the definitions that 
\[  \sup_{y\in Y} \underline{\mdim}\left(\pi^{-1}(y), T\right) \leq \sup_{y\in Y} \overline{\mdim}\left(\pi^{-1}(y), T\right)
     \leq \mdim(\pi, T). \]
So it is enough to prove 
\[  \mdim(\pi, T) \leq \sup_{y\in Y} \underline{\mdim}\left(\pi^{-1}(y), T\right).  \]
Take any positive number $a$ with
\[  \sup_{y\in Y} \underline{\mdim}\left(\pi^{-1}(y), T\right) < a. \]
We will show $\mdim(\pi ,T) \leq  a$.

Let $\varepsilon$ be any positive number.
For every $y\in Y$ there exists $N_y>0$ satisfying $\widim_\varepsilon\left(\pi^{-1}(y), d_{N_y}\right) < a N_y$.

\begin{claim}  \label{claim: extension of epsilon-embedding}
We can find an open neighborhood $U_y$ of $\pi^{-1}(y)$ satisfying 
$\widim_\varepsilon\left(U_y, d_{N_y}\right) < a N_y$. 
\end{claim}
\begin{proof}
Take an $\varepsilon$-embedding $f:\pi^{-1}(y)\to K$ with a simplicial complex $K$ of dimension smaller than $a N_y$.
Since a simplicial complex is ANR (absolute neighborhood retract), 
we can find an open neighborhood $U_y$ of $\pi^{-1}(y)$
and extend $f$ to a continuous map $f:U_y\to K$. 
If we choose $U_y$ sufficiently small, then $f:U_y\to K$ is also an $\varepsilon$-embedding
and we have $\widim_\varepsilon\left(U_y, d_{N_y}\right) < a N_y$.
\end{proof}

There exists an open neighborhood $V_y$ of $y$ satisfying $\pi^{-1}(V_y) \subset U_y$.  
Then we have $\widim_\varepsilon\left(\pi^{-1}(V_y), d_{N_y}\right) < a N_y$.

Since $Y$ is compact, we can find an open cover $Y=V_1\cup V_2 \cup \dots \cup V_m$ and natural numbers $N_1, N_2, \dots, N_m$ satisfying 
$\widim_\varepsilon\left(\pi^{-1}(V_i), d_{N_i}\right) < a N_i$ for all $1\leq i \leq m$.
Set 
\[  \bar{N} := \max_{1\leq i\leq m} N_i. \]

\begin{claim}
For every $y\in Y$ and natural number $N$ we have 
\[  \widim_\varepsilon\left(\pi^{-1}(y), d_N\right) < a (N+\bar{N}). \]
\end{claim}

\begin{proof}
We choose a sequence $i_1, i_2, \dots, i_k$ such that 
\[   y\in V_{i_1}, \quad S^{N_{i_1}}y\in V_{i_2}, \quad S^{N_{i_1}+N_{i_2}}y\in V_{i_3},  \quad \dots \quad ,
     S^{N_{i_1}+N_{i_2}+\dots + N_{i_{k-1}}} y \in V_{i_k}, \]
\[  N_{i_1}+N_{i_2}+\dots+N_{i_{k-1}} < N \leq N_{i_1}+N_{i_2}+\dots + N_{i_k}. \]     
Consider a map 
\[  f: \left(\pi^{-1}(y), d_N\right) \to \left(\pi^{-1}(V_{i_1}), d_{N_{i_1}}\right) \times  \left(\pi^{-1}(V_{i_2}), d_{N_{i_2}}\right) \times \dots
        \times  \left(\pi^{-1}(V_{i_k}), d_{N_{i_k}}\right), \]
defined by $f(x) = \left(x, T^{N_{i_1}}x, T^{N_{i_1}+N_{i_2}}x, \dots, T^{N_{i_1}+N_{i_2}+\dots+N_{i_{k-1}}}x\right)$.
This is a distance non-decreasing continuous map. 
Hence by Lemmas \ref{lemma: product and widim} and \ref{lemma: distance non-decreasing map and widim} we have 
\begin{align*}
   \widim_\varepsilon\left(\pi^{-1}(y), d_N\right) &\leq \sum_{j=1}^k \widim_\varepsilon\left(\pi^{-1}(V_{i_j}), d_{N_{i_j}}\right)  \\
    &< a \sum_{j=1}^k N_{i_j} < a (N+\bar{N}).
\end{align*}    
\end{proof}

Thus we have 
\[  \sup_{y\in Y} \widim_\varepsilon\left(\pi^{-1}(y), d_N\right) \leq a (N+\bar{N}). \]
Hence 
\[  \lim_{N\to \infty} \frac{\sup_{y\in Y} \widim_\varepsilon\left(\pi^{-1}(y), d_N\right)}{N} \leq a. \]
Since $\varepsilon>0$ is arbitrary, we have 
\[  \mdim\left(\pi, T\right) 
    = \lim_{\varepsilon \to 0} \left( \lim_{N\to \infty} \frac{\sup_{y\in Y} \widim_\varepsilon\left(\pi^{-1}(y), d_N\right)}{N} \right)
    \leq a. \]
This proves Proposition \ref{prop: relative mean dimension}.

\section{Proof of Theorem \ref{theorem: Hurewicz does not hold in general}}  
\label{section: proof of Theorem Hurewicz does not hold}

\subsection{Variations of Gromov’s lemma}  \label{subsection: variations of Gromov’s lemma}

Gromov proved the following statement in \cite[p.107, $(\mathrm{H^{\prime\prime}_1})$ Example]{Gromov width}.
This shows that a direct analogue of the Hurewicz theorem does not hold for $\varepsilon$-width dimension.

\begin{lemma}[Gromov 1988]  \label{lemma: Gromov}
Let $(X,d)$ be a $(2n+1)$-dimensional compact Riemannian manifold.
For any positive number $\varepsilon$ there exists a smooth map 
$f:X\to [0,1]$ such that for every $t\in [0,1]$ we have 
\[  \widim_\varepsilon \left(f^{-1}(t), d\right) \leq n.  \]
\end{lemma}

The purpose of this subsection is to develop some variations of this lemma.
The main results are Lemma \ref{lemma: multidimensional Gromov lemma} and Corollary \ref{cor: multidimensional Gromov lemma} below.
We do not explain the proof of Lemma \ref{lemma: Gromov} itself because we will provide a full proof of a more detailed version below.
Our argument follows Gromov’s idea.

We need to prepare some basic terminologies of simplicial complex.
(Recall that we always assume that a simplicial complex has only finitely many simplices.)
For a simplicial complex $K$ we denote by $V(K)$ the set of 
vertices of $K$.

\begin{definition}
Let $K$ be a simplicial complex and $L\subset K$ a subcomplex.
$L$ is said to be a \textbf{full subcomplex} of $K$ if for every simplex $\Delta\subset K$ with $V(\Delta) \subset L$ we have 
$\Delta\subset L$.
\end{definition}

There is a one-to-one correspondence between full subcomplexes of $K$ and a subset of $V(K)$.
For a subset $A\subset V(K)$ we denote by $K(A)$ the (unique) full subcomplex of $K$ satisfying $V\left(K(A)\right) = A$.
If $A = \emptyset$ (empty set), then $K(A) = \emptyset$.

For a vertex $v\in K$ we define the \textbf{star} $\st(v)$ by
\[  \st(v) = \bigcup\left\{\mathrm{Int}(\Delta)\, \middle| \, \text{$\Delta$ is a simplex of $K$ with $v\in \Delta$}\right\}, \]
where $\mathrm{Int}(\Delta)$ is the interior of $\Delta$ defined by (letting $V(\Delta) = \{p_0, p_1, \dots, p_n\}$)
\[  \mathrm{Int}(\Delta) = \left\{\sum_{i=0}^n t_i p_i \, \middle| \, t_0+t_1+\dots+t_n=1, \>  t_i>0 \> (\forall 0\leq i \leq n)\right\}. \] 
When $\Delta = \{v\}$ (just one point) then its interior is $\{v\}$. 
So in particular $\st(v)$ contains $v$.
Indeed it is easy to check that $\st(v)$ is an open neighborhood of $v$.
So the stars $\st(v)$ $(v\in V(K))$ form an open covering of $K$.

Let $m$ be a natural number.
We denote the standard basis of $\mathbb{R}^m$ by $e_1, e_2,\dots, e_m$.
We define the standard $(m-1)$-dimensional simplex $\mathbf{\Delta}^{m-1} \subset \mathbb{R}^m$ by 
\[  \mathbf{\Delta}^{m-1} = \left\{\sum_{i=1}^m t_i e_i\, \middle|\, t_1+t_2+\dots+t_m = 1, \> t_i \geq 0 \> (\forall 1\leq i\leq m) \right\}. \]

The following lemma is our first variation of Gromov’s lemma.

\begin{lemma} \label{lemma: preliminary Gromov lemma}
Let $K$ be a simplicial complex with a  metric $d$. Let $\varepsilon$ be a positive number satisfying
\[  \max_{v\in V(K)} \diam \,\st(v) < \varepsilon. \]
Let $m$ be a natural number.
Suppose we are given a partition 
\[  V(K) = A_1\cup A_2 \cup\cdots \cup A_m \quad (\text{disjoint union}). \]
Then there exists a simplicial map $f:K\to \mathbf{\Delta}^{m-1}$ such that for every $p\in \mathbf{\Delta}^{m-1}$ we have 
\[   \widim_\varepsilon\left(f^{-1}(p), d\right) \leq \max_{1\leq i \leq m} \dim K(A_i). \]
Here $K(A_i)$ is the full subcomplex of $K$ corresponding to $A_i$.
\end{lemma}

\begin{proof}
We define a simplicial map $f:K\to \mathbf{\Delta}^{m-1}$ by the condition $f(A_i) = \{e_i\}$ and extending it linearly.
Namely, if $x\in K$ has the form 
\[  x = \sum_{i=1}^m \left(\sum_{u\in A_i} x_u u\right), \]
where $x_u$ are nonnegative numbers with $\sum_{i=1}^m \sum_{u\in A_i} x_u = 1$, we define 
\[  f(x) = \sum_{i=1}^m \left(\sum_{u\in A_i} x_u\right) e_i. \]

Let $p = \sum_{i=1}^m t_i e_i \in \mathbf{\Delta}^{m-1}$.
We assume, say, $t_1>0$.
(Other cases can be treated similarly.)
If $x = \sum_{i=1}^m \left(\sum_{u\in A_i} x_u u\right) \in f^{-1}(p)$ then 
\[  f(x) = \sum_{i=1}^m \left(\sum_{u\in A_i} x_u\right) e_i = \sum_{i=1}^m t_i e_i, \]
and hence
\[  \sum_{u\in A_1} x_u = t_1 >0. \]
We define a map $g: f^{-1}(p) \to K(A_1)$ by 
\[  g(x) = \frac{\sum_{u\in A_1} x_u u}{\sum_{u\in A_1} x_u}. \]
Every fiber of $g$ is contained in some star $\st(v)$ of $v\in A_1$.
So its diameter is smaller than $\varepsilon$.
Hence $g$ is an $\varepsilon$-embedding.
Therefore 
\[  \widim_\varepsilon\left(f^{-1}(p), d\right) \leq \dim K(A_1). \]
\end{proof}

We need to recall the terminologies on barycentric subdivision.
Let $K$ be a simplicial complex. 
For each simplex $\Delta\subset K$ we denote the barycenter of $\Delta$ by $\bc(\Delta)$.
(Namely, if $V(\Delta) = \{p_0, p_1, \dots, p_n\}$ then $\bc(\Delta) = (p_0+p_1+\dots+p_n)/(n+1)$.)

\begin{figure}[htbp]
\begin{center}
\includegraphics[scale=0.7]{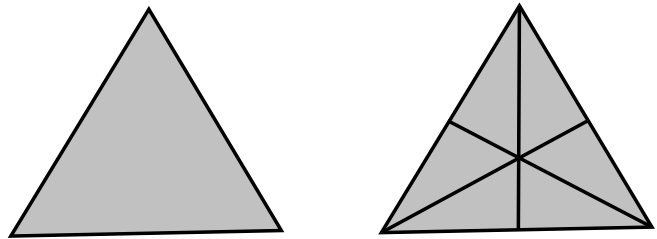} 
\caption{The 2-simplex (left) and its barycentric subdivision (right).}  \label{figure: barycentric subdivision}
\end{center}
\end{figure}

We define the \textbf{barycentric subdivision} $K^\prime$ of $K$ by the following two conditions.
\begin{itemize}
   \item $V\left(K^\prime\right) := \left\{\bc(\Delta)\middle|\, \Delta\subset K: \text{simplex}\right\}$.
   \item If $\Delta_0 \subset \Delta_1 \subset \dots \subset \Delta_n$ is a flag of mutually distinct simplices of $K$ then 
            $\bc(\Delta_0), \bc(\Delta_1), \dots, \bc(\Delta_n)$ form an $n$-simplex in $K^\prime$.
\end{itemize}
As a topological space, $K^\prime$ is naturally identified with $K$.
See Figure \ref{figure: barycentric subdivision}.

The next lemma is our second variation of Gromov’s lemma.

\begin{lemma} \label{lemma: multidimensional Gromov lemma}
Let $K$ be a simplicial complex with a metric $d$. Let $\varepsilon$ be a positive number and $m$ a natural number.
After subdividing $K$ sufficiently fine, we can find a simplicial map 
$f:K\to \mathbf{\Delta}^{m-1}$ such that for every point $p\in \mathbf{\Delta}^{m-1}$
\[  \widim_\varepsilon\left(f^{-1}(p), d\right) \leq \frac{\dim K}{m}. \]
\end{lemma}

\begin{proof}
By subdividing $K$ sufficiently fine, we can assume that 
\[ \max_{v\in V(K)} \diam\,  \st(v) < \varepsilon. \]
Let $K^\prime$ be the barycentric subdivision of $K$.
Then for each simplex $\Delta\subset K$ we have a vertex $\bc(\Delta)$ of $K^\prime$.
We define a partition $V\left(K^\prime\right) = A_1\cup A_2\cup \dots \cup A_m$ (disjoint union) by 
\begin{align*}
  & A_1 = \left\{\bc(\Delta)\middle|\,\text{$\Delta$ is a simplex of $K$ with } \dim \Delta \leq  \frac{\dim K}{m} \right\}, \\
  & A_i = \left\{\bc(\Delta)\middle|\,\text{$\Delta$ is a simplex of $K$ with } \frac{(i-1)\dim K}{m}< \dim \Delta \leq
     \frac{i \dim K}{m} \right\}, \quad (2\leq i\leq m). 
\end{align*}     
Consider the full subcomplexes $K^\prime(A_i)$ $(1\leq i\leq m)$ corresponding to $A_i$.
Every simplex of $K^\prime(A_1)$ corresponds to a flag $\Delta_0\subset \Delta_1\subset \dots \subset \Delta_n$ of 
distinct simplices of $K$ of dimension $\leq \frac{\dim K}{m}$.
The length of such a flag is at most $\frac{\dim K}{m} +1$.
So we have 
\[ \dim K^\prime(A_1) \leq  \frac{\dim K}{m}. \]
Similarly  
\[  \dim K^\prime(A_i) < \frac{\dim K}{m}, \quad (2\leq i \leq m). \]
Figure \ref{figure: Gromov lemma} shows the case that $K$ is two dimensional and $m=2$.

\begin{figure}[htbp]
\begin{center}
\includegraphics[scale=0.5]{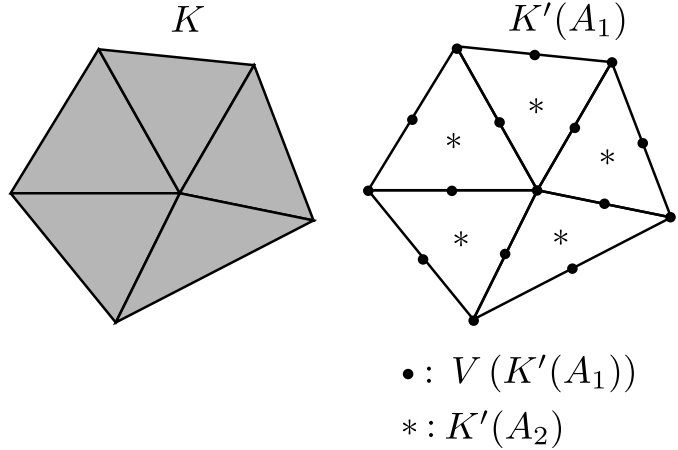} 
\caption{This shows the construction in the case that $K$ is two dimensional and $m=2$.
The left is a two-dimensional simplicial complex $K$.
The right is $K^\prime(A_1)$ and $K^\prime(A_2)$. 
$K^\prime(A_1)$ is the one-dimensional skeleton of $K$ 
(more precisely, $K^\prime(A_1)$ is the barycentric subdivision of the one-dimensional skeleton of $K$).
The vertices of $K^\prime(A_1)$ is depicted as dots.
$K^\prime(A_2)$ consists of five points depicted as $*$, which are the barycenters of two-dimensional simplices of $K$.
So $K^\prime(A_1)$ is one dimensional and $K^\prime(A_2)$ is zero dimensional.}  \label{figure: Gromov lemma}
\end{center}
\end{figure}

Then by Lemma \ref{lemma: preliminary Gromov lemma} we can find 
a simplicial map $f:K^\prime\to \mathbf{\Delta}^{m-1}$ such that for every point $p\in \mathbf{\Delta}^{m-1}$ we have 
\[  \widim_\varepsilon\left(f^{-1}(p), d\right) \leq \max_{1\leq i\leq m} \dim K^\prime(A_i) \leq \frac{\dim K}{m}. \]
\end{proof}

The following corollary is a crucial ingredient of the proof of Theorem \ref{theorem: Hurewicz does not hold in general}
in the next subsection.

\begin{corollary} \label{cor: multidimensional Gromov lemma}
Let $K$ be a simplicial complex with a metric $d$.
For any positive number $\varepsilon$ and any natural number $m$ there exists a continuous map 
$F:K\to [0,1]^{m-1}$ such that for every point $p\in [0,1]^{m-1}$ we have 
\[  \widim_\varepsilon\left(F^{-1}(p), d\right) \leq \frac{\dim K}{m}. \]
\end{corollary}

\begin{proof}
The $(m-1)$ dimensional cube $[0,1]^{m-1}$ is homeomorphic to the $(m-1)$ dimensional simplex $\mathbf{\Delta}^{m-1}$.
So the statement follows from Lemma \ref{lemma: multidimensional Gromov lemma}.
\end{proof}

\begin{remark}
The factor of $\frac{1}{m}$ in the statements of Lemma \ref{lemma: multidimensional Gromov lemma} and 
Corollary \ref{cor: multidimensional Gromov lemma} is optimal, according to 
\cite[p.107, Corollaries $(\mathrm{H^\prime_1})$]{Gromov width}.
It says that if $f:X \to [0,1]^{m-1}$ is a continuous map from a compact metric space $(X, d)$ then for any $\varepsilon>0$
\[  \sup_{p\in [0,1]^{m-1}} \widim_\varepsilon\left(f^{-1}(p), d\right) \geq \frac{\widim_\varepsilon(X, d) -m+1}{m}. \]
\end{remark}

\begin{remark} \label{remark: Gromov lemma is surprising}
Lemma \ref{lemma: multidimensional Gromov lemma} and Corollary \ref{cor: multidimensional Gromov lemma}
show a rather counter-intuitive phenomena.
The following example illustrates its significance:
Let $K$ be a simplicial complex of 
\[ \dim K = 1000,000 \quad  \text{(one million)}. \] 
Let $m=1001$.
From Corollary \ref{cor: multidimensional Gromov lemma}, for any positive number $\varepsilon$, 
there exists a continuous map $F:K\to [0,1]^{1000}$ such that for all $p\in [0,1]^{1000}$ 
\[  \widim_\varepsilon\left(F^{-1}(p), d\right)  < 1000. \]
Hence every fiber of $F$ looks like a space whose dimension is smaller than 1000 (up to distortion bounded by $\varepsilon$).
The range of $F$ is also 1000 dimensional. However the total space $K$ has dimension one million!
\end{remark}

\subsection{Proof of Theorem \ref{theorem: Hurewicz does not hold in general}} 
\label{subsection: Proof of Theorem Hurewicz does not hold}

The purpose of this subsection is to prove Theorem \ref{theorem: Hurewicz does not hold in general}.

Consider the two-sided infinite product of copies of the unit interval $[0,1]$:
\[ [0, 1]^\mathbb{Z} = \cdots \times [0,1]\times [0,1]\times [0,1]\times \cdots. \]
We define a metric $\rho$ on $[0,1]^\mathbb{Z}$ by 
\[ \rho\left((x_n)_{n\in \mathbb{Z}}, (y_n)_{n\in \mathbb{Z}}\right) = \sum_{n\in \mathbb{Z}} 2^{-|n|} |x_n-y_n|. \]
We define the shift map $\sigma:[0,1]^\mathbb{Z}\to [0,1]^\mathbb{Z}$ by 
\[  \sigma\left((x_n)_{n\in \mathbb{Z}}\right) = (x_{n+1})_{n\in \mathbb{Z}}. \]
The pair $\left([0,1]^\mathbb{Z}, \sigma\right)$ is a dynamical system.

Let $(Z, R)$ be a zero dimensional free minimal dynamical system.
Here “\textit{zero dimensional}” means that $Z$ is totally disconnected,
“\textit{free}” means that it has no periodic point, and
“\textit{minimal}” means that every orbit is dense in $Z$.
(Indeed, the following argument works well without the minimality assumption.
But we assume it for simplicity of the explanation.)
It is a standard fact that such a dynamical system exists.
For example, we can construct it as a subshift of $\left(\{0,1\}^\mathbb{Z},\mathrm{shift}\right)$.
We take a metric $\rho^\prime$ on $Z$.

We define a dynamical system $(X, T)$ as the product of $\left([0,1]^\mathbb{Z}, \sigma\right)$ and $(Z,R)$:
\[  (X, T) :=  \left([0,1]^\mathbb{Z}\times Z, \sigma\times R\right). \]
We define a metric $d$ on $X = [0,1]^\mathbb{Z}\times Z$ by 
\[  d\left((x, z), (x^\prime, z^\prime)\right) = \max\left(\rho(x,x^\prime), \rho^\prime(z,z^\prime)\right). \]
The mean dimension of $(X, T)$ is one: 
\[  \mdim(X, T) = 1. \]
Throughout this subsection we fix $(X, T) =  \left([0,1]^\mathbb{Z}\times Z, \sigma\times R\right)$
and construct factor maps from this $(X, T)$.

The next proposition is a preliminary version of Theorem \ref{theorem: Hurewicz does not hold in general}.
\begin{proposition}   \label{prop: Hurewicz does not hold in a fixed scale}
For any positive numbers $\varepsilon$ and $\delta$ there exist a dynamical system $(Y, S)$ and a factor map 
$\pi:(X, T) \to (Y, S)$ such that 
\[  \mdim(Y, S) < \delta, \quad \lim_{N\to \infty} \frac{\sup_{y\in Y}\widim_\varepsilon\left(\pi^{-1}(y), d_N\right)}{N} < \delta. \]
\end{proposition}

\begin{proof}
We prepare ingredients of the construction:
\begin{itemize}
   \item  Fix a natural number $m$ with $\frac{1}{m} < \delta$.
   \item  Fix a natural number $L > m$ with $\frac{m}{L} < \frac{\delta}{2}$.
   \item  Since $(Z, R)$ is free and zero dimensional, we can take a non-empty clopen set 
            $U\subset Z$ satisfying $U\cap R^{-n}U = \emptyset$ for all $1\leq n \leq L$. 
             Here “\textit{clopen}” means that $U$ is both closed and open. Since $(Z, R)$ is minimal, 
             there exists a natural number $L^\prime>L$ satisfying 
            \[  Z = \bigcup_{n=1}^{L^\prime-1} R^{-n} U. \] 
   \item   For $z\in Z$ we set $E(z) = \{n\in \mathbb{Z}\mid R^n z\in U\}$. For any two distinct points $a, b\in E(z)$ we have 
              $|b-a| > L$. Moreover for any $a\in E(z)$ there exists $b\in E(z)$ satisfying $a+L<b<a+L^\prime$.
   \item   For each natural number $n$, we apply Corollary \ref{cor: multidimensional Gromov lemma} to the $n$-dimensional cube $[0,1]^n$.  
              Then we find a continuous map 
              $F_n:[0,1]^n\to [0,1]^{m-1}$ such that for every point $p\in [0,1]^{m-1}$ we have 
              \[  \widim_{\varepsilon/4}\left(F_n^{-1}(p), \norm{\cdot}_\infty\right)  \leq \frac{n}{m}  \left(< \delta n\right). \]
              Here a metric on $[0,1]^n$ is given by the $\ell^\infty$-norm $\norm{x}_\infty = \max_{1\leq i\leq n} |x_i|$.
   \item   For $n\geq m$ we define $G_n:[0,1]^n \to [0,1]^n$ by 
              \[   G_n(x) = \left(F_n(x), \underbrace{0, 0, \dots,0, 0}_{n-m+1}\right). \]           
              Notice that if $n$ is large then the proportion of non-zero entries of $G_n(x)$ is very small.
              For every point $p\in [0,1]^n$ we have 
              \[  \widim_{\varepsilon/4}\left(G_n^{-1}(p), \norm{\cdot}_\infty\right)  \leq \frac{n}{m}  \left(< \delta n\right). \]
\end{itemize}   

We will construct an equivariant continuous map $f:(X, T)\to \left([0,1]^\mathbb{Z},\sigma\right)$
by using the above data.
For a point $x= (x_n)_{n\in \mathbb{Z}} \in [0,1]^\mathbb{Z}$ and integers $a<b$, we denote 
\[  x|_{[a, b)} := (x_a, x_{a+1}, x_{a+2}, \dots, x_{b-1}).  \]
Let $(x, z) \in [0,1]^\mathbb{Z}\times Z = X$.
Take any point $a\in E(z)$ and set $b := \min \left(E(z) \cap (a, \infty)\right)$.
We have $L<b-a<L^\prime$.
We define
\[  f(x, z)|_{[a, b)} := G_{b-a}\left(x|_{[a,b)}\right) \in [0,1]^{b-a}. \]
We consider this for every $a\in E(z)$. Then we have defined $f(x, z)\in [0,1]^\mathbb{Z}$.
The map $f$ is equivariant and continuous. (The continuity follows from the clopenness of $U$.)
We define $\pi:X\to [0,1]^\mathbb{Z}\times Z$ by 
\[ \pi(x, z) = \left(f(x,z), z\right). \]
We set $Y = \pi(X)$, which becomes a dynamical system under the map $S:= \sigma\times R$.
We will show that a factor map $\pi:(X, T) \to (Y, S)$ satisfies the statement.

First we estimate the mean dimension of $(Y, S)$.
Let $(x,z)\in [0,1]^\mathbb{Z}\times Z$ and let $N$ be a natural number.
Denote 
\[  [0, N)\cap E(z) = \{a_1<a_2<a_3<\dots<a_k\}. \]
Since $a_{i+1}-a_i >L$, we have $k -1< \frac{N}{L}$.
Set 
\[ a_0 := \max\left(E(z) \cap (-\infty, 0)\right), \quad a_{k+1} := \min\left(E(z)\cap [N, \infty)\right). \]
Then we have 
\[  [0,N) \subset \bigcup_{i=0}^k [a_i, a_{i+1}) = [a_0, a_{k+1}). \]
For $0\leq i \leq k$
\[   f(x, z)|_{[a_i, a_{i+1})} = G_{a_{i+1}-a_i}\left(x|_{[a_i, a_{i+1})}\right). \]
The number of its non-zero entries is at most $m-1$.  Therefore
\begin{align*}
    \text{The number of non-zero entries of $f(x,z)|_{[0,N)}$} &\leq (k+1) (m-1)  \\
    & < \left(\frac{N}{L} + 2\right) (m-1)  \\
    & < \frac{\delta N}{2} + 2m  \quad \text{by $\frac{m}{L} < \frac{\delta}{2}$}. 
\end{align*}
Denote by $\Pi_N: [0,1]^\mathbb{Z}\to [0,1]^N$ the projection to the $0, 1, 2, \dots, (N-1)$-th coordinates.
Then $\Pi_N\left(f(X)\right)$ is contained in 
\[  \left\{\left(y_0, y_1, \dots, y_{N-1}\right)\, \middle|\, 
     y_n = 0 \text{ except for at most $\frac{\delta N}{2} + 2m$ entries}\right\}. \]
whose dimension is at most $\frac{\delta N}{2} + 2m$.
This implies that 
\[ \mdim\left(f(X), \sigma\right) \leq \lim_{N\to \infty} \left(\frac{\delta N}{2} + 2m\right)/N = \frac{\delta}{2} < \delta. \]
Since $Z$ is zero dimensional, we have 
\[  \mdim(Y, S) = \mdim\left(f(X),\sigma\right)  <\delta. \]

Next we study the fibers of $\pi$.
Fix $M>0$ with $\sum_{|n|\geq M} 2^{-|n|} < \frac{\varepsilon}{2}$. 
Let $y = (p, z)\in Y$ with $p\in [0,1]^\mathbb{Z}$ and $z\in Z$.
Let $N$ be a natural number.
Denote 
\[  (-M, N+M)\cap E(z) = \{a_1<a_2<\dots<a_k\},  \]
and set 
\[  a_0 := \max\left(E(z)\cap (-\infty, -M]\right), \quad a_{k+1} := \min\left(E(z)\cap [N+M, \infty)\right). \]
Then we have 
\[  (-M, N+M) \subset \bigcup_{i=0}^k [a_i, a_{i+1}) = [a_0, a_{k+1}). \]
Since $a_{i+1}-a_i < L^\prime$, we have $a_{k+1}-a_0 < N+2M+2L^\prime$.
For any $0\leq i\leq k$ and $(x, z)\in \pi^{-1}(y)$ we have 
\[   f(x,z)|_{[a_i, a_{i+1})} = G_{a_{i+1}-a_i}\left(x|_{[a_i, a_{i+1})}\right) = p|_{[a_i,a_{i+1})}.   \]

We define a projection $\Pi_{[a_0, a_{k+1})}: [0,1]^\mathbb{Z} \times Z \to [0,1]^{a_{k+1}-a_0}$ by 
\[  \Pi_{[a_0, a_{k+1})}\left((x_n)_{n\in \mathbb{Z}}, z\right) = (x_n)_{n\in [a_0, a_{k+1})}. \]
Then 
\[   \Pi_{[a_0, a_{k+1})}\left(\pi^{-1}(y)\right) \subset
    G_{a_1-a_0}^{-1}\left(p|_{[a_0, a_1)}\right) \times G_{a_2-a_1}^{-1}\left(p|_{[a_1, a_2)}\right)\times \cdots \times G_{a_{k+1}-a_k}^{-1}\left(p|_{[a_k, a_{k+1})}\right).   \]
It follows that 
\begin{align*}
    \widim_\varepsilon\left(\pi^{-1}(y), d_N\right) & 
    \leq \widim_{\varepsilon/4}\left(\Pi_{[a_0, a_{k+1})}\left(\pi^{-1}(y)\right), \norm{\cdot}_\infty\right) \quad
     \text{by the definition of the metric $d$} \\
    & \leq  \sum_{i=0}^k \widim_{\varepsilon/2} \left(G_{a_{i+1}-a_i}^{-1}\left(p|_{[a_i, a_{i+1})}\right), \norm{\cdot}_\infty\right) \\
    & \leq \sum_{i=0}^k  \frac{a_{i+1}-a_i}{m} \\
    & = \frac{a_{k+1}-a_0}{m} < \frac{N+2M+2L^\prime}{m}.
\end{align*}
This holds for every $y\in Y$. So we get 
\[ \sup_{y\in Y}  \widim_\varepsilon\left(\pi^{-1}(y), d_N\right)  < \frac{N+2M+2L^\prime}{m}. \]
Notice that $M$ and $L^\prime$ are independent of $N$.
Thus we conclude 
\[  \lim_{N\to \infty} \frac{\sup_{y\in Y}  \widim_\varepsilon\left(\pi^{-1}(y), d_N\right)}{N} \leq \frac{1}{m} < \delta. \]
\end{proof}

Now we are ready to prove Theorem \ref{theorem: Hurewicz does not hold in general}.

\begin{theorem} \label{theorem: Hurewicz does not hold version two}
For any positive number $\delta$ there exists a dynamical system $(Y, S)$ and a factor map 
$\pi:(X, T)\to (Y, S)$ such that 
\[   \mdim(Y, S) < \delta, \quad \mdim(\pi, T) = 0. \]
\end{theorem}

Recall that the dynamical system $(X, T) = \left([0,1]^\mathbb{Z}\times Z, \sigma\times R\right)$ has mean dimension one.
So Theorem \ref{theorem: Hurewicz does not hold in general} follows from this theorem.

\begin{proof}
For each natural number $n$
we apply Proposition \ref{prop: Hurewicz does not hold in a fixed scale} to $\varepsilon = \frac{1}{n}$ and $\frac{\delta}{2^n}$.
Then we find a factor map $\pi_n:(X, T) \to (Y_n, S_n)$ such that 
\[ \mdim(Y_n, S_n) < \frac{\delta}{2^n}, \quad 
    \lim_{N\to \infty} \frac{\sup_{y\in Y_n} \widim_{1/n}\left(\pi_n^{-1}(y), d_N\right)}{N} < \frac{\delta}{2^n}. \]
Define 
\[  \pi:= \pi_1\times \pi_2\times \pi_3\times\cdots : X\to Y_1\times Y_2\times Y_3\times \cdots, 
\quad x\mapsto (\pi_1(x), \pi_2(x), \pi_3(x), \dots). \]
Set $Y = \pi(X) \subset Y_1\times Y_2\times Y_3\times \cdots$
with a map $S:= S_1\times S_2\times S_3\times \cdots$.
The pair $(Y, S)$ is a dynamical system.
We show that the factor map $\pi:(X, T)\to (Y, S)$ satisfies the statement.

The mean dimension of $(Y, S)$ is bounded by 
\[ \mdim(Y, S) \leq \sum_{n=1}^\infty \mdim(Y_n, S_n) < \sum_{n=1}^\infty \frac{\delta}{2^n} = \delta. \]

Let $y = (y_1, y_2, y_3, \dots)\in Y$ $(y_n\in Y_n)$.
We have 
\[  \pi^{-1}(y)  = \bigcap_{n=1}^\infty \pi_n^{-1}(y_n).  \]
Then for any natural numbers $n$ and $N$
\[  \widim_{1/n} \left(\pi^{-1}(y),d_N\right) \leq \widim_{1/n}\left(\pi_n^{-1}(y_n), d_N\right). \]  
Hence 
\[  \lim_{N\to \infty} \frac{\sup_{y\in Y} \widim_{1/n}\left(\pi^{-1}(y), d_N\right)}{N} 
     \leq \lim_{N\to \infty} \frac{\sup_{y_n\in Y_n}  \widim_{1/n}\left(\pi_n^{-1}(y_n), d_N\right)}{N} < \frac{\delta}{2^n}. \]
Letting $n\to \infty$, we get $\mdim(\pi, T) =0$.
\end{proof}

\begin{remark}  \label{remark: raising mean dimension}
Let $\delta$ be a positive number. For each natural number $n$ we apply Theorem \ref{theorem: Hurewicz does not hold version two}
to $\frac{\delta}{2^n}$. Then there exists a factor map $\pi_n: (X, T) \to (Y_n, S_n)$ such that 
\[  \mdim(Y_n, S_n) < \frac{\delta}{2^n}, \quad \mdim(\pi_n, T) = 0.   \]
We set
\[  (X^\prime, T^\prime) := (X, T) \times (X, T) \times (X,T) \times \cdots, \quad 
     (Y^\prime, S^\prime) := (Y_1, S_1) \times (Y_2, S_2) \times (Y_3, S_3) \times \cdots. \]
We define a factor map $\pi^\prime: (X^\prime, T^\prime) \to  (Y^\prime, S^\prime)$ by 
\[  \pi^\prime(x_1, x_2, x_3, \dots) = (\pi_1(x_1), \pi_2(x_2), \pi_3(x_3), \dots). \]
Then it is easy to check that 
\[  \mdim(X^\prime, T^\prime) = \infty, \quad \mdim(Y^\prime, S^\prime) < \delta, \quad \mdim(\pi^\prime, T^\prime) = 0. \]
See Remark \ref{remark: raising mean dimension of (X,T)} in \S \ref{subsection: main results}.
\end{remark}

\begin{remark} \label{remark: universality}
The dynamical system $(X^\prime, T^\prime)$ constructed in the above remark is a rather “universal” one.
It has the form 
\[  (X^\prime, T^\prime) = \left(\left([0,1]^\mathbb{N}\right)^\mathbb{Z}, \mathrm{shift}\right) \times (Z^\prime, S^\prime) \]
where $(Z^\prime, S^\prime)$ is a zero dimensional dynamical system given by 
\[  (Z^\prime, S^\prime) = (Z, S) \times (Z, S) \times (Z, S) \times \cdots. \]
It is easy to see that \textit{every} dynamical system embeds in $\left(\left([0,1]^\mathbb{N}\right)^\mathbb{Z}, \mathrm{shift}\right)$.
Hence, given an arbitrary dynamical system $(\mathcal{X}, \mathcal{T})$, we can embed 
$\left(\mathcal{X}, \mathcal{T}\right) \times (Z^\prime, S^\prime)$ in $(X^\prime, T^\prime)$.
Therefore, for any positive number $\delta$, there exists 
a factor map 
\[  \Pi: \left(\mathcal{X} \times Z^\prime, \mathcal{T}\times S^\prime\right) \to \left(\mathcal{Y}, \mathcal{S}\right) \]
satisfying 
\[  \mdim\left(\mathcal{Y}, \mathcal{S}\right) < \delta, \quad \mdim\left(\Pi, \mathcal{T}\right) = 0. \]
Notice that 
\[  \mdim \left(\mathcal{X} \times Z^\prime, \mathcal{T}\times S^\prime\right) = \mdim\left(\mathcal{X},\mathcal{T}\right) \]
can be an arbitrary nonnegative number.
Thus we can say that the construction of this section shows an universal phenomena.
\end{remark}

\section{Proof of Theorem \ref{theorem: Hurewicz holds if the base is zero mean dimensional}}
\label{section: proof of Theorem Hurewicz holds}

\subsection{Preliminaries on relative mean dimension}  \label{subsection: preliminaries on relative mean dimension}

Here we prepare some simple facts on the relative mean dimension.

\begin{lemma}  \label{lemma: relative mean dimension of product}
 Let $\pi_i: (X_i, T_i)\to (Y_i, S_i)$ $(i=1, 2)$ be two factor maps between dynamical systems.
 We consider their product:
 \[  \pi_1\times \pi_2: (X_1\times X_2, T_1\times T_2)\to (Y_1\times Y_2, S_1\times S_2). \]
 For this factor map we have 
\begin{align*}
    &    \mdim\left(\pi_1\times \pi_2, T_1\times T_2\right) \leq  \mdim(\pi_1,T_1) + \mdim(\pi_2, T_2),  \\
    &   \mdim\left(\pi_1\times \pi_2, T_1\times T_2\right)  \geq \max\left(\mdim(\pi_1,T_1), \mdim(\pi_2,T_2)\right).
\end{align*}    
\end{lemma}

\begin{proof}
Let $d$ and $d^\prime$ be metrics on $X_1$ and $X_2$ respectively.
The product space $X\times X^\prime$ has a metric $d\times d^\prime$. (See \S \ref{subsection: preparations on widim}.)

Let $(y_1,y_2)\in Y_1\times Y_2$. For every natural number $N$ we have 
\[ \left((\pi_1\times \pi_2)^{-1}(y_1,y_2), (d\times d^\prime)_N\right) = \left(\pi_1^{-1}(y_1), d_N\right) \times \left(\pi_2^{-1}(y_2), d^\prime_N\right). \]
By Lemma \ref{lemma: product and widim}, for $\varepsilon>0$
\[  \widim_\varepsilon\left((\pi_1\times \pi_2)^{-1}(y_1,y_2), (d\times d^\prime)_N\right) \leq 
     \widim_\varepsilon\left(\pi_1^{-1}(y_1), d_N\right)   + \widim_\varepsilon\left(\pi_2^{-1}(y_2), d^\prime_N\right). \]
Thus we have the first inequality.

Fix $p\in \pi_2^{-1}(y_2)$. The map 
\[  \left(\pi_1^{-1}(y_1), d_N\right) \to  \left((\pi_1\times \pi_2)^{-1}(y_1,y_2), (d\times d^\prime)_N\right), \quad 
      x \mapsto (x, p), \]
is an isometric embedding. Hence 
\[ \widim_\varepsilon\left(\pi_1^{-1}(y_1), d_N\right) 
    \leq \widim_\varepsilon\left((\pi_1\times \pi_2)^{-1}(y_1,y_2), (d\times d^\prime)_N\right). \]
Similarly 
\[  \widim_\varepsilon\left(\pi_2^{-1}(y_2), d_N\right) 
    \leq \widim_\varepsilon\left((\pi_1\times \pi_2)^{-1}(y_1,y_2), (d\times d^\prime)_N\right). \]
Then we get the second inequality.
\end{proof}

\begin{corollary}  \label{corollary: extension trick}
 Let $\pi:(X, T)\to (Y, S)$ be a factor map between dynamical systems, and let $(Z, R)$ a dynamical system.
 We consider 
 \[  \pi \times \mathrm{Id}: (X\times Z, T\times R) \to (Y\times Z, S\times R), \quad 
       (x, z) \mapsto (\pi(x), z). \]
Then we have 
\[  \mdim\left(\pi\times \mathrm{Id}, T\times R \right) = \mdim(\pi, T). \]
\end{corollary}

\begin{proof}
The trivial factor map 
\[  \mathrm{Id}: (Z, R) \to (Z, R) \]
has zero relative mean dimension. So the corollary follows from Lemma \ref{lemma: relative mean dimension of product}.
\end{proof}

\subsection{Small boundary property}  \label{subsection: small boundary property}

Here we review the theory of \textit{small boundary property} introduced by Lindenstrauss--Weiss \cite{Lindenstrauss--Weiss}.
Let $(X, T)$ be a dynamical system.
For a subset $A\subset X$ we define the \textbf{orbit capacity} of $A$ by 
\[ \ocap(A) = \lim_{N\to \infty} \frac{1}{N} \sup_{x\in X} \sum_{n=0}^{N-1} 1_A(T^n x). \]
The quantity $\sup_{x\in X} \sum_{n=0}^{N-1} 1_A(T^n x)$ is sub-additive in $N$.
So we have 
\begin{equation}  \label{eq: orbit capacity}
    \ocap(A) = \inf_{N\geq 1} \frac{1}{N} \sup_{x\in X} \sum_{n=0}^{N-1} 1_A(T^n x). 
\end{equation}    
It is easy to see that $\ocap(A\cup B) \leq \ocap(A) + \ocap(B)$.
In particular, if $\ocap(A) = \ocap(B)=0$ then $\ocap(A\cup B) =0$.

The next lemma was proved in \cite[Lemma 6.3]{Lindenstrauss}.

\begin{lemma} \label{lemma: ocap and perturbation}
Let $E\subset X$ be a closed subset. For any positive number $\delta$ there exists an open neighborhood $U$ of $E$ satisfying 
\[   \ocap(U) < \ocap(E)+\varepsilon. \]
\end{lemma}

\begin{proof}
This follows from the formula (\ref{eq: orbit capacity}).
\end{proof}

We say that a dynamical system $(X, T)$ has the \textbf{small boundary property} if for every point $x\in X$ and for every open 
neighborhood $U$ of $x$ there exists an open set $V$ such that 
\[  x\in V \subset U, \quad \ocap(\partial V) = 0. \]
Here $\partial V$ is the boundary of $V$, namely $\partial V := \overline{V}\setminus V$.
The small boundary property is a dynamical version of totally disconnectedness.

Lindenstrauss--Weiss \cite[Theorem 5.4]{Lindenstrauss--Weiss} proved that if a dynamical system has the 
small boundary property then its mean dimension is zero.
Lindenstrauss \cite[Theorem 6.2]{Lindenstrauss} proved a partial converse of this statement as follows. 
This will be crucial in the next subsection.

\begin{theorem}[Lindenstrauss 1999]  \label{theorem: small boundary property}
If $(X, T)$ is an extension of a free minimal system and $\mdim(X, T) = 0$ then 
$(X, T)$ has the small boundary property.
\end{theorem}

Here the assumption that “$(X, T)$ is an extension of a free minimal system” means that 
there exists a factor map $\pi:(X, T) \to (Y, S)$ such that $(Y,S)$ is a free minimal dynamical system.

\begin{lemma} \label{lemma: open cover}
Let $(Y, S)$ be a dynamical system having the small boundary property.
Let $\delta$ be a positive number.
For any open covering $Y = V_1\cup V_2\cup \dots \cup V_m$ there exist compact subsets $E_i\subset V_i$ $(1\leq i\leq m)$
such that 
\begin{itemize}
   \item $E_i \cap E_j = \emptyset$ for $i\neq j$.
   \item $\ocap\left(Y\setminus (E_1\cup \dots\cup E_m)\right) < \delta$.
\end{itemize}
\end{lemma}

\begin{proof}
By the small boundary property there exist open sets $W_i$ $(1\leq i\leq m)$ such that 
$\overline{W_i}\subset V_i$, $\ocap(\partial W_i) =0$ and $Y = W_1\cup W_2\cup\dots \cup W_m$.
By Lemma \ref{lemma: ocap and perturbation} there exist open sets $U_i \supset \partial W_i$ with 
$\ocap(U_i) < \frac{\delta}{m}$ for $1\leq i\leq m$.
Then $\overline{W_i}\cup U_i = W_i \cup U_i$ is an open set.
We set 
\begin{align*}
   & E_1 := \overline{W_1}, \quad E_2 := \overline{W_2}\setminus \left(\overline{W_1}\cup U_1\right), \quad 
    E_3 := \overline{W_3}\setminus \left(\overline{W_1}\cup U_1 \cup \overline{W_2}\cup U_2\right), \quad \dots, \\
  &  E_m := \overline{W_m}\setminus \left(\overline{W_1}\cup U_1\cup \overline{W_2}\cup U_2\cup \dots \cup \overline{W_{m-1}}\cup U_{m-1}\right).
\end{align*}    
$E_i$ are compact and 
\[ Y  \setminus (E_1\cup \dots\cup E_m)  \subset U_1\cup U_2\cup \dots \cup U_m. \]
Then
\[ \ocap\left(Y\setminus (E_1\cup \dots\cup E_m)\right)  \leq \sum_{i=1}^m \ocap(U_i) < \delta. \]
\end{proof}

\subsection{Proof of Theorem \ref{theorem: Hurewicz holds if the base is zero mean dimensional}}  
\label{subsection: proof of Theorem Hurewicz holds}

In this subsection we first
prove an analogue of the Hurewicz theorem for mean dimension in the case that the base system has the small boundary property.
Next we prove it in the case that the base system has zero mean dimension.

Before stating the proposition, we prepare terminologies on \textit{cone}.
This will be used in the proof.
Let $K$ be a topological space. We define the \textbf{cone} $C(K)$ by 
\[  C(K) = [0,1]\times K/\sim, \]
where $(0, x)\sim (0, y)$ for any $x,y\in K$.
For $0\leq t\leq 1$ and $x\in X$, the equivalence class of $(t,x)$ is denoted by $tx$.
The point $0x$ is called the \textbf{vertex} of the cone $C(K)$ and often denoted by $*$.

When $K$ is a simplicial complex,
the cone $C(K)$ naturally admits a structure of a simplicial complex. Its dimension is $\dim K +1$.

\begin{figure}[htbp]
\begin{center}
\includegraphics[scale=0.7]{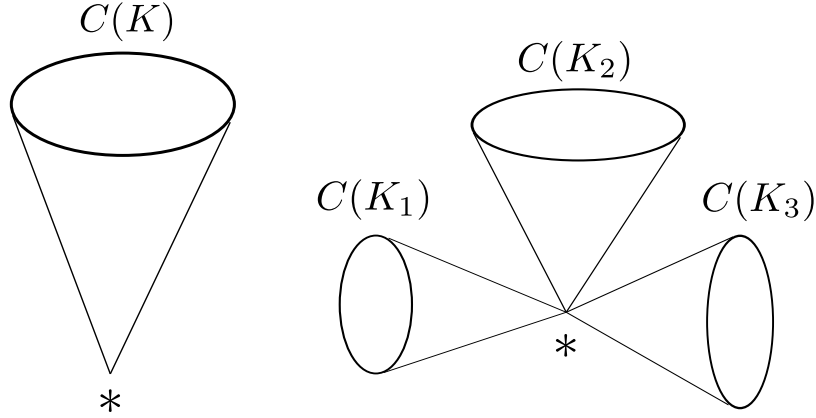} 
\caption{The left is the cone $C(K)$ with the vertex $*$. The right is $C(K_1)\cup_* C(K_2)\cup_* C(K_3)$ with the vertex $*$. 
Three cones are glued at their vertices.}  \label{figure: cone}
\end{center}
\end{figure}

Let $K_1, \dots, K_m$ be topological spaces and let $C(K_1), \dots, C(K_m)$ their cones.
We denote by $*_i$ the vertex of the cone $C(K_i)$.
Let $C(K_1)\cup C(K_2)\cup\dots \cup C(K_m)$ be the disjoint union of $C(K_1), \dots, C(K_m)$.
We define 
\[  C(K_1)\cup_* C(K_2) \cup_*\dots \cup_* C(K_m) = C(K_1)\cup C(K_2)\cup\dots \cup C(K_m)/\sim \]
where $*_i \sim *_j$ for all $i, j$.
Namely we glue $C(K_1), \dots, C(K_m)$ at their vertices, and the resulting space is denoted by
$C(K_1)\cup_* C(K_2) \cup_*\dots \cup_* C(K_m)$. The shared vertex (i.e. the equivalence class of $*_i$) is denoted by $*$.
See Figure \ref{figure: cone}.

When $K_1, \dots, K_m$ are simplicial complexes then 
$C(K_1)\cup_*\dots \cup_* C(K_m)$ is also a simplicial complex and its dimension is the maximum of $\dim K_i + 1$ $(1\leq i\leq m)$.

\begin{proposition}  \label{prop: Hurewicz and small boundary property}
Let $\pi:(X, T)\to (Y, S)$ be a factor map between dynamical systems.
If $(Y,S)$ has the small boundary property then
\[   \mdim(X,T) = \mdim(\pi, T). \]
\end{proposition}

\begin{proof}
$\mdim(\pi, T) \leq \mdim(X,T)$ is obvious from the definition.
Let $a$ be any positive number with $\mdim(\pi, T) < a$.
We will prove $\mdim(X,T)\leq a$.

We take a metric $d$ on $X$.
Let $\varepsilon$ be any positive number.
We can take a natural number $N$ such that 
\[  \sup_{y\in Y} \widim_\varepsilon\left(\pi^{-1}(y), d_N\right) +1 < a N. \]
As in Claim \ref{claim: extension of epsilon-embedding} in \S \ref{subsection: proof of Proposition relative mean dimension},
for each $y\in Y$ there exists an open neighborhood $U_y$ of $\pi^{-1}(y)$ satisfying 
\[  \widim_\varepsilon(U_y, d_N) +1 < a N. \]
There exists an open neighborhood $V_y$ of $y$ with $\pi^{-1}(V_y) \subset U_y$.

Since $Y$ is compact, we can find an open covering $Y = V_1\cup V_2\cup \dots \cup V_m$ such that for all 
$1\leq i\leq m$
\[   \widim_\varepsilon\left(\pi^{-1}(V_i), d_N\right) + 1 < a N.  \]
Let $\delta$ be any positive number.
By applying Lemma \ref{lemma: open cover}
to an open cover $Y = V_1\cup \dots\cup V_m$, we find compact subsets $E_i\subset V_i$ $(1\leq i \leq m)$ such that 
\[  E_i\cap E_j = \emptyset \quad (i\neq j), \quad 
     \ocap\left(Y\setminus (E_1\cup \dots \cup E_m)\right) < \delta. \]
We take open sets $W_i$ $(1\leq i \leq m)$ satisfying 
\[  E_i \subset W_i \subset V_i,  \quad  W_i \cap W_j = \emptyset \quad (i\neq j). \]
We also take a continuous function $\rho:Y \to [0,1]$ such that $\rho = 1$ on $E_1\cup E_2\cup \dots \cup E_m$ and 
$\supp \, \rho \subset W_1\cup W_2\cup \dots \cup W_m$.
See Figure \ref{figure: cutoff}.

\begin{figure}[htbp]
\begin{center}
\includegraphics[scale=0.7]{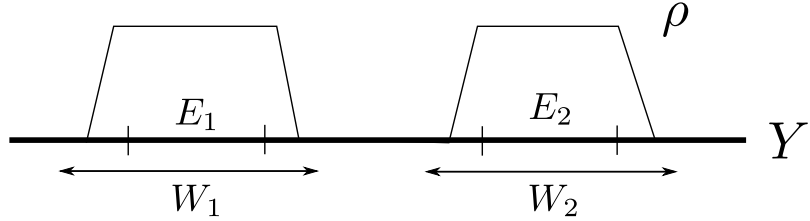} 
\caption{The continuous function $\rho$.}  \label{figure: cutoff}
\end{center}
\end{figure}

For each $1\leq i\leq m$ we take an $\varepsilon$-embedding 
$f_i:\left(\pi^{-1}(V_i), d_N\right) \to K_i$ such that $K_i$ is a simplicial complex of dimension smaller than 
$aN-1$.
Let $C(K_i)$ be the cone over $K_i$. 
We set
\[  K^\prime := C(K_1)\cup_* C(K_2) \cup_* \dots \cup_* C(K_m).  \]
This is a simplicial complex of dimension smaller than $aN$.
We define a continuous map $f^\prime:X\to K^\prime$ as follows:
If $x\in \pi^{-1}(W_i)$ $(1\leq i\leq m)$ then we set 
\[  f^\prime(x) = \rho\left(\pi(x)\right) f_i(x) \in C(K_i) \subset K^\prime.     \]
If $x\not\in \pi^{-1}(W_1)\cup \pi^{-1}(W_2)\cup\cdots\cup \pi^{-1}(W_m)$ then we set 
$f^\prime(x) = *$.

We take an $\varepsilon$-embedding $g:(X, d_N) \to L$ such that $L$ is a simplicial complex of dimension 
$\widim_\varepsilon(X,d_N) < \infty$.
Let $L^\prime := C(L)$ be the cone over $L$.
We have $\dim L^\prime = \widim_\varepsilon(X, d_N)+1$.
We define a continuous map $g^\prime:X\to L^\prime$ by 
\[  g^\prime(x) = \left(1-\rho(\pi(x))\right) g(x).  \]

For $n\geq 1$ we define $F_n:X\to \left(K^\prime\times L^\prime\right)^n$ by 
\[  F_n(x) = \left(f^\prime(x), g^\prime(x), f^\prime(T^N x), g^\prime(T^N x), \dots, f^\prime(T^{(n-1)N} x), g^\prime(T^{(n-1)N} x)\right). \]

\begin{claim}  \label{claim: F_n is epsilon embedding}
$F_n$ is an $\varepsilon$-embedding with respect to the metric $d_{nN}$.
\end{claim}

\begin{proof}
It is enough to prove that 
\[  f^\prime\times g^\prime:X\to K^\prime\times L^\prime, \quad x\mapsto (f^\prime(x), g^\prime(x)) \]
is an $\varepsilon$-embedding with respect to the metric $d_N$.

Suppose $\left(f^\prime(x), g^\prime(x)\right) = \left(f^\prime(x^\prime), g^\prime(x^\prime)\right)$ for $x, x^\prime\in X$.
Then we have $\rho\left(\pi(x)\right) = \rho\left(\pi(x^\prime)\right)$. If this common value is zero then 
$g(x) = g(x^\prime)$ and hence we have $d_N(x, x^\prime) < \varepsilon$.
If the value is positive then there exists $1\leq i\leq m$ such that $x, x^\prime\in \pi^{-1}(W_i)$ and $f_i(x)= f_i(x^\prime)$.
Then $d_N(x,x^\prime) < \varepsilon$.
\end{proof}

\begin{claim}  \label{claim: number of nonzero g}
If $n$ is sufficiently large then for all $x\in X$
\[  \left|\{0\leq k <n \mid g^\prime(T^{kN}x) \neq *\}\right| < \delta n N. \]
Therefore the image of $F_n$ is contained in a simplicial complex of dimension smaller than 
$n \dim K^\prime + \delta n N \dim L^\prime$.
\end{claim}

\begin{proof}
Let $x\in X$ and set $y=\pi(x)$.
The condition $g^\prime\left(T^{kN}x\right) \neq *$ is equivalent to $\rho\left(S^{kN}y\right) <1$.
The latter condition implies $S^{kN}y \in Y\setminus (E_1\cup \dots \cup E_m)$.
Then
\begin{align*}
     \left|\{0\leq k <n \mid g^\prime(T^{kN}x) \neq *\}\right|  & \leq 
     \sum_{k=0}^{n-1} 1_{Y\setminus (E_1\cup \dots \cup E_m)}(S^{kN} y)  \\
     & \leq \sum_{k=0}^{nN-1}1_{Y\setminus (E_1\cup \dots \cup E_m)}(S^{k} y)  \\
     & < \delta nN  \quad \text{by $\ocap \left(Y\setminus (E_1\cup \dots \cup E_m)\right) < \delta$}
\end{align*}
for all sufficiently large $n$ (uniformly in $x \in X$).
\end{proof}

By Claims \ref{claim: F_n is epsilon embedding} and \ref{claim: number of nonzero g}, for sufficiently large $n$
\begin{align*}
    \widim_\varepsilon\left(X, d_{nN}\right)  & < n \dim K^\prime + \delta n N \dim L^\prime \\
    & < anN + \delta n N \left(\widim_\varepsilon(X, d_N) + 1\right).
\end{align*}
Then 
\[ \lim_{n\to \infty} \frac{\widim_\varepsilon\left(X, d_{nN}\right)}{nN} \leq 
    a +  \delta \left(\widim_\varepsilon(X, d_N) + 1\right). \]
Here $\delta$ is independent of $\varepsilon, N$. So we can let $\delta\to 0$ and get 
\[  \lim_{n\to \infty} \frac{\widim_\varepsilon(X, d_n)}{n} \leq a. \]
Letting $\varepsilon\to 0$, we conclude
\[  \mdim\left(X, T\right) \leq a. \]
\end{proof}

Now we are ready to prove Theorem \ref{theorem: Hurewicz holds if the base is zero mean dimensional}.
We write the statement again.

\begin{theorem}[$=$ Theorem \ref{theorem: Hurewicz holds if the base is zero mean dimensional}]
Let $\pi:(X, T)\to (Y,S)$ be a factor map between dynamical systems.
If $\mdim(Y,S) =0$ then
\[   \mdim(X,T) = \mdim(\pi, T). \]
\end{theorem}

\begin{proof}
Let $(Z, R)$ be a zero dimensional free minimal dynamical system. 
We consider 
\begin{equation}  \label{eq: product factor map trick} 
     \pi\times \mathrm{Id}: (X\times Z, T\times R) \to (Y\times Z, S\times R). 
\end{equation}     
The system $(Y\times Z, S\times R)$ has a free minimal factor $(Z,R)$ and its mean dimension is zero:
\[  \mdim(Y\times Z, S\times R) = \mdim(Y,S) =0 \quad    \text{since $\dim Z=0$}. \]
By the Lindenstrauss Theorem (Theorem \ref{theorem: small boundary property}), the system $(Y\times Z, S\times R)$ has the small boundary property.

Now we can apply Proposition \ref{prop: Hurewicz and small boundary property} to the above factor map (\ref{eq: product factor map trick}) and get
\[  \mdim(X\times Z, T\times R) = \mdim\left(\pi\times \mathrm{Id}, T \times R\right). \]
We have 
\[  \mdim(X\times Z, T\times R) = \mdim(X, T) \quad \text{since $\dim Z=0$}, \]
\[  \mdim\left(\pi\times \mathrm{Id}, T \times R\right)  = \mdim(\pi, T) \quad \text{by Corollary \ref{corollary: extension trick}}. \]
Therefore we conclude $\mdim(X, T) = \mdim(\pi, T)$.
\end{proof}

\end{document}